\documentclass{elsart3-1}

\usepackage{amssymb}

\usepackage[english,francais]{babel}

\usepackage{amsmath}
\usepackage{mathrsfs}
\usepackage[all]{xy}
\usepackage[dvips]{graphicx}
\usepackage{url}

\newtheorem{theorem}{Theorem}[section]
\newtheorem{lemma}[theorem]{Lemma}
\newtheorem{e-proposition}[theorem]{Proposition}

\newtheorem{e-definition}[theorem]{Definition\rm}
\newtheorem{remark}{\it Remark\/}

\renewcommand{\theequation}{\arabic{equation}}
\newcommand{\of}{\mathcal{O}}

\newcommand{\pn}{\mathbb{P}}
\newcommand{\rk}{\mathrm{rk}\,}

\def\og{\leavevmode\raise.3ex\hbox{$\scriptscriptstyle\langle\!\langle$~}}
\def\fg{\leavevmode\raise.3ex\hbox{~$\!\scriptscriptstyle\,\rangle\!\rangle$}}

\journal{the Acad\'emie des sciences}

\begin{document}

\begin{frontmatter}

\selectlanguage{english}
\title{About the stability of the tangent bundle restricted to a curve}

\selectlanguage{english}
\author[authorlabel1]{Chiara Camere},
\ead{chiara.camere@unice.fr}

\address[authorlabel1]{Laboratoire J.-A. Dieudonn\'e U.M.R. no 6621 du C.N.R.S.
Universit\'e de Nice-Sophia Antipolis
Parc Valrose
06108 Nice}

\begin{abstract}
\selectlanguage{english}
Let $C$ be a smooth projective curve of genus $g\geq 2$ and let $L$ be a line bundle on $C$ generated by its global sections. The morphism $ \phi _L:C\longrightarrow\pn (H^0(L))\simeq\pn ^r $ is well-defined and $\phi _L^*T_{\pn ^r}$ is the restriction to $C$ of the tangent bundle of $\pn ^r$. Sharpening a theorem by Paranjape, we show that if $\deg L\geq 2g-c(C)$ then $\phi _L^*T_{\pn ^r}$ is semi-stable, specifying when it is also stable. We then prove the existence on many curves of a line bundle $L$ of degree $2g-c(C)-1$ such that $\phi _L^*T_{\pn ^r}$ is not semi-stable. Finally, we completely characterize the (semi-)stability of $\phi _L^*T_{\pn ^r}$ when $C$ is hyperelliptic.

\vskip 0.5\baselineskip

\selectlanguage{francais}
\noindent{\bf R\'esum\'e} \vskip 0.5\baselineskip \noindent
{\bf Sur la stabilit\'e du fibr\'e tangent restreint \`a une courbe. }
Soit $ L $ un fibr\'e en droites engendr\'e par ses sections globales sur une courbe projective lisse $ C $ de genre $ g\geq 2 $. Le fibr\'e $L$ d\'efinit $ \phi _L:C\longrightarrow\pn (H^0(L))\simeq\pn ^r $ et $\phi _L^*T_{\pn ^r}$ est la restriction \`a la courbe $C$ du fibr\'e tangent de $\pn ^r$. En pr\'ecisant un th\'eor\`eme d\^u \`a Paranjape, on montre que si $\deg L\geq 2g-c(C)$ alors $\phi _L^*T_{\pn ^r}$ est semi-stable, en disant quand il est aussi stable. De plus, on montre l'existence sur plusieurs courbes d'un fibr\'e en droites $L$ de degr\'e $2g-c(C)-1$ tel que $\phi _L^*T_{\pn ^r}$ ne soit pas semi-stable. Enfin, on caract\'erise compl\`etement la stabilit\'e de $\phi _L^*T_{\pn ^r}$ si $C$ est hyperelliptique.

\end{abstract}
\end{frontmatter}

\selectlanguage{english}
\section{Introduction}
\label{}
Let $C$ be a smooth projective curve of genus $g\geq 2$ and let $L$ be a line bundle on $C$ generated by its global sections. Let $M_L$ be the vector bundle defined by the exact sequence
\begin{equation}\label{eq:syz}
\xymatrix{
0\ar@{->}[r]&M_L\ar@{->}[r]&H^0(C,L)\otimes\mathcal{O}_C \ar@{->}[r]^-{e_L}&L\ar@{->}[r]&0}
\end{equation} 
where $e_L$ is the evaluation map. We denote by $E_L$ the dual bundle of $M_L$: it has degree $\deg L$ and rank $h^0(C,L)-1$.
Let us briefly recall the geometric interpretation of these bundles: since $L$ is generated by its global sections, the morphism $ \phi _L:C\longrightarrow\pn (H^0(L))\simeq\pn ^r $ is well-defined and we have $L=\phi _L^*\of _{\pn ^r}(1) $; thus, from the dual sequence of  (\ref{eq:syz}) and from the well-known Euler exact sequence
\begin{equation} \label{s:eu}
 \xymatrix{
0\ar@{->}[r]&\of _{\pn ^r}\ar@{->}[r]&H^0(C,L)^*\otimes\of _{\pn ^r}(1)\ar@{->}[r]&T _{\pn ^r}\ar@{->}[r]&0} 
\end{equation}
it follows that $ E_L=\phi _L^*T_{\pn ^r}\otimes L^* $ and the stability of $E_L $ is equivalent to the stability of $\phi _L^*T_{\pn ^r}$.

We recall the definition of the Clifford index of a curve.
\begin{e-definition}
The Clifford index of a line bundle $L$ on $C$ is $c(L)=\deg L-2(h^0(C,L)-1).$

The Clifford index of a divisor $D$ on $C$ is the Clifford index of the associated line bundle $\of _C(D)$, i.e. $ c(D)=c(\of _C(D))=\deg D-2\dim |D| .$

The Clifford index of the curve $C$ is $ c(C)=\min \{c(L)/ h^0(C,L)\geq 2,\  h^1(C,L)\geq 2\} .$
\end{e-definition}

Clifford's theorem states that $c(C)\geq 0$, with equality if and only if $C$ is hyperelliptic; moreover, for any divisor $D$ on $C$, $ c(D)=c(K-D) $.
\begin{remark}\label{r:cl}
 By the Riemann-Roch theorem, $c(L)=2g-\deg L-2h^1(C,L)$ for any line bundle $L$.
\end{remark}

In \cite{P}, by using the properties of this invariant, Paranjape proves the following 
\begin{e-proposition}
Let $C$ be a smooth projective curve of genus $g\geq 2$ and let $L$ be a line bundle on $C$ generated by its global sections. If $c(C)\geq c(L)$ then $E_L$ is semi-stable. If $h^1(C,L)=1$ and $ c(C)> 0 $ or $c(C)> c(L)$ then $E_L$ is also stable.
\end{e-proposition}

By completing his proof we show the following
\begin{theorem}\label{p:pa}
Let $C$ be a smooth projective curve of genus $g\geq 2$ and let $L$ be a line bundle on $C$ generated by its global sections such that $\deg L\geq 2g-c(C) $. Then:
\begin{enumerate}
 \item $ E_L $ is semi-stable;
 \item $ E_L $ is stable except when $\deg L= 2g$ and either $C$ is hyperelliptic or $L\cong K(p+q)$ with $p,q\in C$.
\end{enumerate}
\end{theorem}

If $C$ is a smooth projective $d-$gonal curve of genus $g\geq 2$ with Clifford index $c(C)=d-2<\frac{g-2}{2}$, we then prove the existence of a line bundle $L$ of degree $2g-c(C)-1$ such that $E_L$ is not semi-stable. Moreover, a theorem by Schneider (see \cite{Sc}) states that on a general smooth curve $E_L$ is always semi-stable: our proof also shows that one cannot replace semi-stable by stable in this statement.

Finally, we completely characterize the (semi-)stability of $E_L$ when $C$ is hyperelliptic.

\section{Proof of Theorem \ref{p:pa}}
We first need a lemma, shown by Paranjape in \cite{P}.

\begin{lemma}\label{p:p32}
Let $F$ be a vector bundle on $C$ generated by its global sections and such that $ H^0(C,F^*)=0 $; then 
$ \deg F\geq \rk F+g-h^1(C,\det F) $ and equality holds if and only if $ F=E_L $, where $ L=\det F $.
Moreover, if $ h^1(C,\det F)\geq 2 $ then $ \deg F\geq 2\rk F+c(C) $ and if equality holds then $ F=E_L $.
\end{lemma}

The canonical bundle $K$ is generated by its global sections and there is an exact sequence
\[ 
\xymatrix{
0\ar@{->}[r]&K^*\ar@{->}[r]&H^0(C,K)^*\otimes \of _C\ar@{->}[r]&E_K\ar@{->}[r]&0
} 
 \]
thus in cohomology we have
\begin{equation}\label{s:ka}
 \xymatrix{
0\ar@{->}[r]&H^0(K^*)\ar@{->}[r]& H^0(K)\!^{^*}\!\!\otimes\! H^0(\of _{_C})\ar@{->}[r]& H^0(E _{_K})\ar@{->}[r]&H^1(K^*)\ar@{->}[r]^{\varphi\ \ \ \ \ \ }&H^0(K)\!^{^*}\!\!\otimes\! H^1(\of _{_C})\ar@{->}[r]&\cdots
} 
\end{equation}
The map $ \varphi $ is the dual map of $m:H^0(K)\otimes H^0(K)\rightarrow H^0(K^2)$, so it is injective by Noether's theorem (see \cite{ACGH}, Chap.III); moreover, $ H^0(C, K^*)=0 $. As a consequence $ H^0(C,E_K)\simeq H^0(C,K)^*=H^1(C,\of _C) $ and $ h^0(C,E_K)=g $.

Now we have all the tools necessary to prove Theorem \ref{p:pa}.

\textbf{Proof of Theorem \ref{p:pa}.} By Remark \ref{r:cl}, if $\deg L\geq 2g-c(C) $ a fortiori $c(C)\geq c(L) $. By definition, $ \deg E_L= c(L)+2\rk E_L$ and $ h^0(C,L)=\rk E_L+1 $, hence it follows by the Riemann-Roch theorem that $\deg E_L= \rk E_L+g-h^1(C,L).$

Let $F$ be a quotient bundle of $E_L$; then $F$ satisfies the hypothesis of Lemma \ref{p:p32}, because it is spanned by its global sections since $E_L$ is and $ H^0(C,F^*)\subset H^0(C,E_L^*)=0 $. 

Therefore, if $ h^1(C,\det F)\geq 2 $ we have $ \deg F\geq 2\rk F+c(C) $; then \[ \mu (F)-\mu (E_L)\!\geq\!\frac{c(C)}{\rk F}-\frac{c(L)}{\rk E_L}=\frac{\rk E_L\!\cdot\! c(C)-\rk F\!\cdot\! c(L)}{\rk F\cdot \rk E_L}=\frac{(\rk E_L-\rk F)\!\cdot\! c(C)+\rk F\!\cdot\!(c(C)-c(L))}{\rk F\cdot\rk E_L}\!\geq\! 0 \] since $ \rk E_L> \rk F >0$ and $ c(C)\geq c(L)  $. Moreover, the inequality is strict if $ c(C)> 0 $ or if $C$ is hyperelliptic and $\deg L\geq 2g+1$, because $L$ is non-special and $c(L)<0$.

If $ h^1(C,\det F)< 2 $ we still have $  \deg F\geq \rk F+g-h^1(C,\det F)  $, hence
\[ \mu (F)-\mu (E_L)\!\geq\!\frac{g\!-\!h^1(\det F)}{\rk F}-\frac{g\!-\!h^1(L)}{\rk E_L}\!=\!\frac{\left[ g\!-\!h^1(\det F)\right]\!\!\cdot\! (\rk E_L\! -\!\rk F)\!+\!\rk F\!\cdot\! \left[ h^1(L)\!-\! h^1(\det F)\right] }{\rk F\cdot\rk E_L}\!>\! 0\]
provided that $ h^1(C,L)\geq h^1(C,\det F)$, since $g-h^1(C,\det F)>0$ follows from the hypothesis that $ h^1(C,\det F)< 2 $ and $g\geq 2$.

The only case remaining is $ 0=h^1(C,L)<h^1(C,\det F)=1.$ We have $ \deg F=\deg(\det F)\leq 2g-2 $, otherwise we should have $ h^1(C,\det F)=0 $; then, a fortiori, we have $ \rk F\leq g-1 $. It then follows from the previous inequalities that 
\begin{equation}\label{e:st}
 \mu (F)-\mu (E_L)\!\geq\! 
\frac{(g-1)(\rk E_L-\rk F)\!-\!\rk F}{\rk F\cdot\rk E_L}\!\geq\!\frac{(g-1)\!\cdot\!(\rk E_L-\rk F-1)}{\rk F\cdot\rk E_L}\!\geq\! 0
\end{equation}
Thus we have shown that we always have $\mu (F)-\mu (E_L)\geq 0$, i.e. $E_L$ is semi-stable. In order to gain the stability of $E_L$, we still need to prove that $\mu (F)-\mu (E_L)> 0$ when $ 0=h^1(C,L)<h^1(C,\det F)=1$.

Suppose that $ \mu (E_L)=\mu (F)$; by (\ref{e:st}), we then have $ (g-1)\!\cdot\!\rk E_L-g\!\cdot\!\rk F=0. $ Since $ g\geq 2 $, it follows that $ (g-1)|\rk F\leq g-1 $, i.e. $ \rk F=g-1 $, and $ \rk E_L=g $; hence $ \deg E_L=g+\rk E_L=2g $ and $ \mu (E_L)=2 $. Therefore, if $\deg L\neq 2g$ we cannot have $ \mu (E_L)=\mu (F)$ and $E_L$ is stable.

If $\deg L= 2g$ then $E_L$ is stable provided that $c(C)>0$ and $L\ncong K(p+q)$ with $p,q \in C$.

Indeed, since $ \deg F=\rk F\!\cdot\! \mu (F)=2g-2 $ and $ h^1(C,\det F)=1 $, we have $ \det F \cong K$. As a consequence we have $\rk F+g-h^1(C,\det F)=2g-2=\deg F,$ so $ F=E_K $ by Lemma \ref{p:p32}. On the other hand, $F$ is a quotient of $E_L$, so there is an exact sequence
\begin{equation}\label{s:w}
\xymatrix{
0\ar@{->}[r]& W\ar@{->}[r]&E_L\ar@{->}[r]&F\ar@{->}[r]&0
} 
\end{equation} 
where $W$ is a sub-bundle of $E_L$ of degree 2 and rank 1. The associated exact sequence of cohomology then is
\[
\xymatrix{
0\!\ar@{->}[r]&\! H^0(C,\!W)\!\ar@{->}[r]&\! H^0(C,E_L)\!\ar@{->}[r]^{\varphi}&\! H^0(C,E_K)\!\ar@{->}[r]&\! H^1(C,\!W)\!\ar@{->}[r]&\!\cdots
}\]
From the exact sequence of cohomology associated to the dual sequence of (\ref{eq:syz}) we see that $ h^0(C,E_L)\geq g+1 $ and $ h^0(C,E_K)=g $ since $c(C)>0$; hence $ \varphi $ cannot be injective, i.e. $ H^0(C,W)\neq 0 $. Thus $W\cong \of _C(p+q) $ with $ p,q\in C $. Furthermore, it follows from (\ref{s:w}) that $$ L=\det E_L=\det W\otimes\det F=W\otimes K =K(p+q),$$ which concludes the proof of Theorem \ref{p:pa} since this is not possible under our hypothesis.\qed

\section{Some line bundles of degree $2g-c(C)-1$ with non semi-stable $E_L$}

Theorem \ref{p:pa} is the best possible result that one can obtain if looking for properties of all curves.

\begin{e-proposition}\label{p:cep}
Let $C$ be a smooth projective d-gonal curve of genus $g\geq 2$ such that the Clifford index is $c(C)=d-2<\frac{g-2}{2}$; there exists a line bundle $L$ of degree $ \deg L=2g-c(C)-1$ on $C$ generated by its global sections and non-special such that $E_L$ is not semi-stable.
\end{e-proposition}
 
\textbf{Proof.} By the hypothesis, $\mathfrak{g}^1 _d$ computes the Clifford index. We put $N=\of _C (K-\mathfrak{g}^1 _d)$: it is a line bundle of degree $2g-c(C)-4$ and by the Riemann-Roch theorem $h^0(N)=g-c(C)-1.$ Moreover $N$ is spanned by its global sections: assume that there exists $q\in C$ such that $h^0(N(-q))=h^0(N)$, or equivalently $h^1(N(-q))=h^1(N)+1$; then, by Serre's duality, we have $h^0(\mathfrak{g}^1 _d+q)=h^0(\mathfrak{g}^1 _d)+1=3$, i.e. $\mathfrak{g}^1 _d+q=\mathfrak{g}^2 _{d+1}$, and this is not possible because we would have $c(\mathfrak{g}^2 _{d+1})=d-3<c(C)$.

Let $E$ be an effective divisor of degree 3 on $C$; we can choose $E$ in such a way that $L=N\otimes \of _C(E)$ is a line bundle of degree $ \deg L=2g-c(C)-1$, non-special and spanned by its global sections. Indeed, we have $h^1(L)=0$ because $h^1(L)=h^0(\mathfrak{g}^1 _d-E)=0$ for a general effective divisor $E$; moreover $L$ is generated by its global sections if and only if $ h^1(L(-p))=h^1(L)=0$ for any $p\in C$ and if $E$ is a general effective divisor of degree 3 we have $h^1(L(-p))=h^0(\mathfrak{g}^1 _d-E+p)=0$.

Since we have supposed that $E$ is effective, $H^0(L\otimes N^*))\neq 0$, so we have an inclusion $N\hookrightarrow L$. Hence $M_N$ is a sub-bundle of $M_L$, or equivalently $E_N$ is a quotient bundle of $E_L$. Since $\rk E_L=g-c(C)-1$ and $\rk E_N=h^0(N)-1=g-c(C)-2$, we have
\begin{equation}\label{d:dg}
 \mu (E_N)=2+\frac{c(C)}{g-c(C)-2}<\mu (E_L)=2+\frac{c(C)+1}{g-c(C)-1}
\end{equation}
whenever $c(C)<\frac{g-2}{2}$. It then follows that $E_L$ is not semi-stable.\qed

\begin{remark}
If $C$ is a curve of genus $g\geq 2$ with Clifford index $c$, in most cases $C$ is $(c+2)-$gonal: see \cite{ELMS} for further details.
\end{remark}

\begin{remark}
The hypothesis that $c(C)<\frac{g-2}{2}$ leaves out only the case $c(C)=\left[ \frac{g-1}{2}\right]$, i.e. the general one; however, in \cite{Sc} Schneider shows the following
\begin{e-proposition}
 Let $C$ be a general smooth curve of genus $g\geq 3$. If $L$ is a line bundle on $C$ generated by its global sections, then $E_L$ is semi-stable.
\end{e-proposition}

It is worth underlining that one cannot replace semi-stable by stable: if $C$ is a general curve of even genus $g=2n$ we know that 
\begin{equation}\label{e:clmax}
 c(C)=\left[ \frac{g-1}{2}\right]=n-1=\frac{g-2}{2},
\end{equation}
so the proof of Proposition \ref{p:cep} shows that $E_L$ is not stable, since one obtains $\mu (E_N)=\mu (E_L)$.
\end{remark}

\section{The case of hyperelliptic curves}

In the case of hyperelliptic curves we completely characterize the stability of $E_L$.

\begin{e-proposition}
Let $C$ be a smooth projective hyperelliptic curve of genus $g\geq 2$, let $L$ be a line bundle on $C$ generated by its global sections and such that $h^0(C,L)\geq 3$ and let $ H$ be $\of _C(\mathfrak{g}^1 _2) $. Then:
\begin{enumerate}
 \item $E_L$ is stable if and only if $\deg L\geq 2g+1 $;
 \item $E_L$ is semi-stable if and only if $ \deg L\geq 2g $ or there exists an integer $k>0$ such that $L=H^{\otimes k}$.
\end{enumerate}
\end{e-proposition}

\textbf{Proof.} By Theorem \ref{p:pa}, if $\deg L\geq 2g$ then $E_L$ is semi-stable and if $\deg L\geq 2g+1 $ then $E_L$ is stable.

On the other hand $E_L$ is not stable if $\deg L=2g$, in which case $ \mu (E_L)=2 $. Indeed, we show that $ H$ is a quotient bundle of $E_L$ of same slope. We know that there is a surjection $ E_L\twoheadrightarrow H $ if and only if there is an inclusion $ H^*\rightarrowtail M_L $, if and only if $ H^0(C,M_L\otimes H)\neq 0 $. From the exact sequence (\ref{eq:syz}) we get an exact sequence
\begin{equation}\label{s:mlh}
\xymatrix{
0\ar@{->}[r]& H^0(C,M_L\otimes H)\ar@{->}[r]&H^0(C,L)\otimes H^0(C,H)\ar@{->}[r]&H^0(C,L\otimes H)\ar@{->}[r]&\cdots
} 
\end{equation}
We then have $ \dim H^0(C,L)\otimes H^0(C,H)=2g+2>g+3= h^0(C,L\otimes H) $, so $ H^0(C,M_L\otimes H)\neq 0 $.

If $0<\deg L\leq 2g-1$ we always have $c(L)\geq 0$. If $c(L)= 0$ then $E_L$ is semi-stable, as it follows from the proof of Theorem \ref{p:pa}: if $F$ is a quotient bundle of $E_L$, the inequality $\mu (F)-\mu (E_L)\geq 0$ still holds in each case.

Using again the exact sequence (\ref{s:mlh}), since $h^0(C,L)\geq 3$, we have \[\dim H^0(C,L)\otimes H^0(C,H)=2h^0(C,L)>h^0(C,L)+2\geq h^0(C,L\otimes H).\]
Therefore, $ H^0(C,M_L\otimes H)\neq 0 $ and there is a surjection $ E_L\twoheadrightarrow H $; furthermore, $$\mu (E_L)=2+\frac{c(L)}{h^0(C,L)-1} $$
and $\mu (H)=2.$ Thus if $c(L)> 0$ then $\mu (E_L)>\mu (H)$ and $E_L$ is not semi-stable; else, if $c(L)=0$, $\mu (E_L)=\mu (H)$ and $E_L$ is not stable.

The proposition then follows by Clifford's theorem: since $C$ is hyperelliptic and $\deg L>0$, $c(L)=0$ if and only if there exists an integer $k>0$ such that $L=H^{\otimes k}$.\qed

\section*{\textbf{Acknowledgements}}

I would like to thank my supervisor Professor Arnaud Beauville for his patient guidance and for the time he spent reading the various drafts of this article.

\end{document}